# A MATHEMATICAL THEORY OF ORIGAMI CONSTRUCTIONS AND NUMBERS

ROGER C. ALPERIN

## 1. Introduction

About twelve years ago, I learned that paper folding or elementary origami could be used to demonstrate all the Euclidean constructions; the booklet, [9], gives postulates and detailed the methods for high school teachers. Since then, I have noticed a number of papers on origami and variations, [3], [2] and even websites [7]. What are a good set of axioms and what should be constructible all came into focus for me when I saw the article [13] on constructions with conics in the Mathematical Intelligencer.

The constructions described here are for the most part classical, going back to Pythagorus, Euclid, Pappus and concern constructions with ruler, scale, compass, and angle trisections using conics. Klein mentioned the book of Row, [12], while describing geometrical constructions in [10], but went no further with it. Row's book uses paper folding, as he says, 'kindergarten tools', to study geometrical constructions and curve sketching.

We shall describe a set of axioms for paper folding which will be used to describe, in a hierarchial fashion, different subfields of the complex numbers, in the familiar way that ruler and compass constructions are used to build fields. The axioms for the origami constructible points of the complex numbers, starting with the constructible points 0 and 1 are that it is the smallest subset of constructible points obtained from the following axioms:

(1) The line connecting two constructible points is a constructible line.

(2) The point of coincidence of two constructible lines is a constructible point.

(3) The perpendicular bisector of the segment connecting two constructible points is a constructible line.

(4) The line bisecting any given constructed angle can be constructed.

(5) Given a constructed line $l$ and constructed points $P, Q$, then whenever possible, the line through $Q$, which reflects $P$ onto $l$, can be constructed.

(6) Given constructed lines $l, m$ and constructed points $P, Q$, then whenever possible, any line which simultaneously reflects $P$ onto $l$ and $Q$ onto $m$, can be constructed.

These axioms allow constructions of lines, which are performed in origami by folding a piece of paper. The constructed points make up the origami numbers. The points on a constructed line are not necessarily constructible points. The first three axioms, which we call Thalian constructions, do not seem very strong at all, using merely perpendicular bisections, but surprisingly starting with a third non-real point give the structure of a field to the constructed set of points. The fourth axiom, allowing angle bisections, in a sense completes the first level giving the Pythagorean numbers, studied by Hilbert in *Foundations of Geometry* in connection with constructions with a ruler and (unit) scale and their relations to the totally real algebraic numbers. In [1], this idea is developed using a larger set of axioms. The fifth axiom, adds yet more constructions, precisely the Euclidean constructions, not by using a compass, but by adding in the construction of the envelope of tangents of a parabola. This has been discussed by [9], [3] with additional axioms and the use of double folds, but it is classical origami and geometrical constructions. In fact as we shall show, one can eliminate axioms (1) and (4) as a application of the power of using axiom (5). The axioms (2), (3), (5) are all that are needed for geometrical constructions. The sixth axiom, has been discussed before in [3], and in [2] using the mira constructions. This axiom allows the constructions of cube roots, solving the problem of the duplication of the cube, just as the ancients did it, using the intersection of parabolas, [13]. This last axiom admits the construction of the tangents to two parabolas as a new construction. This is strong enough to be used to solve any cubic or fourth order equation using resolvant techniques. The method of the cubic resolvant of a fourth order equation is in fact related to the idea of intersections of conics.

Our main contribution here is to show that with all six axioms we get precisely the field obtained from intersections of conics, the field obtained from the rationals by adjoining arbitrary square roots and cube roots and conjugates. The techniques here are elementary algebraic geometry, [3], the theory of pencils of conics or quadratic forms.

Of course, the standard question, as to which regular polygons can be constructed, is readily answered, [13], [2]; however, Gleason, [4], who develops the theory of the angle trisector, also derives the same conclusion, that the number of sides is $2^a 3^b P_1 P_2 .. P_s$, where the distinct primes $P_i$, if any, are of the form $2^c 3^d + 1$.

A good reference for solving equations in one variable, and its history, are contained in [8]. The classical theory of Euclidean constructions is carried out there, too. Also, many interesting and historically relevant comments to classical constructions of cube roots and trisections of angles with the aid of curves are made in all of the references to this paper, so few historical points will be repeated here.

## 2. Geometrical Axioms and Algebraic Consequences

2.1. **Thalian Constructions.** We have named the first collection of axioms after Thales, the teacher of Pythagorus. Thales founded the first Ionian school of mathematics in Miletus in the 6th century BC. The axioms are described so that they can easily be implemented with paper and folds. The axioms are rather weak so we have to work hard to prove anything. But it all falls into place nicely in the end.

We assume the complex numbers $\mathbb{C}$ are given. All constructions are assumed to take place in $\mathbb{C}$. The requirements for a (Thalian) constructible set of numbers $\Pi = \Pi\{A, B, C\}$ are that it is the smallest subset of constructible points in the plane of the complex numbers, which is closed under the following operations:

(0) The set contains 3 points $A$, $B$, $C$, of the complex numbers, not all on a line.

(1) The line connecting two points in $\Pi$ is constructible.

(2) The point of coincidence of two constructible lines is in $\Pi$.

(3) The perpendicular bisector of the segment connecting two constructible points is constructible.

It is easy to see how one can implement these axioms using a single sheet of paper, three marked points, and then either folding the paper to create a constructible line through marked points, or marking the intersection of two folds, or folding two marked points onto each other to create the perpendicular bisector.

For the Lemma below, and its Corollaries, we assume the constructions are taking place in a given $\Pi$.

**Lemma 2.1.** *Given a point $P$ and a line segment $AB$ then we can construct the parallel segment of the same length and direction as $AB$ beginning or ending at $P$.*

**Proof.** Suppose first that $P$ is not on the line through $AB$. Connect $P$ to $A$ and $B$ and bisect the sides of the constructed triangle, obtaining points $p, a, b$ on the sides opposite the given points.

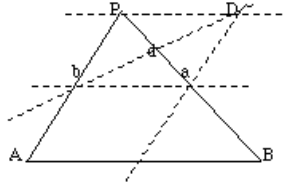

Figure 2.1. Translation

As in Figure 2.1, the line $ab$ is parallel to the base line $BA$; construct the line $pa$. Bisect the segment $Pa$ to obtain point $d$. Construct the line $bd$. The line $bd$ and $pa$ intersect at $D$. The segment $PD$ is parallel to $AB$ and half its length. We can now do a similar construction for the line parallel to $ap$ and passing through $B$ to obtain the point on the line $PD$ of the desired length $AB$.

The construction for the segment ending at $P$ is done similarly. If we are given a point $P$ on the line through $AB$, then we can similarly move the segment to start or end at a point $Q$, not on the line, and then move back to $P$. ∎

**Corollary 2.2.** *We can construct a perpendicular to a line $m$ from a given point $P$.*

**Proof.** Using two constructible points on $m$ construct the perpendicular bisector using axiom 3 and now its parallel through $P$, by Lemma 2.1. ∎

**Corollary 2.3.** *Given a point $P$ and a segment $ABC$ we can construct the segment $ADP$ so that triangle $ABD$ is directly similar to $ACP$.*

**Proof.** Construct the parallel through $B$ of the line through $PC$. If $P$ is on the segment $ABC$ first make the construction as before with an arbitrary point $Q$ not on the segment $ABC$ of the segment $AEQ$.

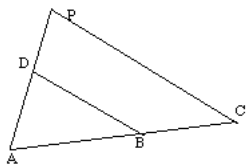

Figure 2.2. Similarity

Now with the segment $AEQ$, and the point $P$, we may find the point $D$ so that $ADP$ is in the same ratio as $AEQ$ which is the same as $ABC$. Of course in this case $D$ is on the segment $ABC$. ∎

**Corollary 2.4.** *Given a point $P$ and a line $l$ we can reflect $P$ across $l$. Given lines $l$ and $m$ we can reflect $m$ across $l$.*

**Proof.** The second construction includes the first.. Take two constructible points on $l$, $A$ and $B$. Construct the perpendicular to $l$ at $A$ and the intersection with $m$, $a$. Similarly with $B$. Translate the segment $Aa$ to $A$, and get new endpoint $C$; do similarly for $B$ to get $D$. This gives two reflected points $C$ and $D$ and then construct the line passing through them. ∎

**Corollary 2.5.** *The subset $\Pi$ is closed under segment addition in $\mathbb{C}$. If we take $A=0$ then $\Pi$ is an abelian group.*

If we assume that the first given point $A$ is 0 and the second given point $B$ is the complex number 1, then we can construct the x-axis and also its perpendicular y-axis. The set $\Pi\{0, 1, C\}$ depends now only on one complex variable $z$, the third point $C$; we shall denote the set as $\Pi = \Pi[z]$. In this case we may construct the x-coordinates and y-coordinates as subsets of the real numbers, $X = X[z], Y = Y[z]$. Since we can project a constructed point to the x-axis or y-axis, $X[z]$ and $iY[z]$ are subsets of $\Pi$. From the abelian group structure on $\Pi$, we have the structure of abelian groups on $X$ and $Y$ as subsets of the real numbers and moreover, $\Pi[z] = X[z] \oplus iY[z]$.

**Corollary 2.6.** *For $z$ a non-real complex number, $\Pi[z]$ is a $\mathbb{Q}$-vector space closed under complex conjugation, with subspaces and $X[z]$ and $iY[z]$.*

**Proof.** Given the point $W$, and a positive integer $n$ we can solve for $U$ so that $nU = W$. First construct $nV$ for some vector not on the same line as $W$. Construct the line through $nV$ and $W$; the parallel through $V$ passes through $U$. The remaining claims follow immediately from previous remarks. ∎

**Lemma 2.7.** *If $t \in Y$ is non-zero then $1/t \in Y$.*

**Proof.** We may asume $t > 0$. We make the construction of a right triangle with legs $t$ along the y-axis and 1 parallel to the x-axis. The hypotenuse begins at the origin and extends to $(1, t)$. Make a second triangle: drop a perpendicular to the y-axis from $(1, t)$ and also construct a perpendicular to the hypotenuse at $(1, t)$ extending to the y-axis. By similar triangles, the length of the leg along the y-axis of the second triangle is $1/t$.

These subsets $X, Y \subset \Pi$ have more structure. With two points $u, x \in X$ and $y \in Y$ we construct $(x, y)$ and get by Corollary 2.3, $(u, uy/x)$ so that $uy/x \in Y$. If $u = 1$ then $y/x \in Y$. Similarly, if $v, y \in Y$, and $x \in X$, we construct $(x, y)$ and by Corollary 2.3 we obtain $(vx/y, v)$ so that $vx/y \in X$. If $x = 1$, then $v/y \in X$. (One can think of these as multiplying real numbers $X$, and imaginary numbers $Y$.) If $z$ is not real then $\Pi[z]$ has a non-zero $y \in Y$.

**Corollary 2.8.** Let $x \in X$, $t, v, y \in Y$:
  i) $vy \in X$, $xy \in Y$ and consequently, $tvy \in Y$;
  ii) $x^2 \in X$, if $Y \neq \{0\}$, hence $X$ is a $\mathbb{Q}$-algebra;
  iii) for non-zero $y \in Y$, then $Y = Xy$;
  iv) for non-zero $x \in X$, if $Y \neq \{0\}$, then $1/x \in X$, and hence $X$ is a field;
  v) if $Y \cap X \neq \{0\}$ then $X = Y$;
  vi) if $X = Y$ then $\Pi[z]$ is the field, $X(i)$.
  vii) for any real number $\mu$, $\mu$ is the slope of a constructed line iff $\mu \in Y$.

**Proof.** The conclusions about $vy$ follow immediately from the remarks above and the following. Let $1 = u, x \in X$, $y \in Y$, then by Lemma 2.7, $1/y \in Y$ and therefore $1/(xy) \in Y$; hence by Lemma 2.7 again $xy \in Y$. Thus $Xy \subset Y$. Again following the remarks above, with $v, y \in Y$, $x \in X$, $xy \in Y$, so $vx/(xy) \in X$. Thus $v \in Xy$, so that $Y \subset Xy$. For ii) construct $(x, y)$. The line from the origin passing through $(x, y)$ meets the line parallel to the x-axis at height $xy$, a constructible in $Y$, at the point $(x^2, xy)$ so $x^2 \in X$. Since $X$ contains $\mathbb{Q}$, and is a $\mathbb{Q}$ vector space, closure under multiplication follows from $2uv = (u + v)^2 - u^2 - v^2$. For iv) let $y \in Y$ be non-zero, $x$ a non-zero element of $X$, then $1/(yx) \in Y$ by i) and Lemma 2.7; also $y^2 \in X$, by i) and so $y/x = y^2/(yx) \in Y$, by i). Therefore $1/x = (1/y)(y/x)$ is in $X$ by i) and Lemma 2.7. Thus, it follows from iv) and ii) that $X$ is a field. For part v) if $y \in Y \cap X$ is non-zero then $1/y \in Y$, so $1 \in Y$, and thus $Y = X$.

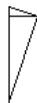

FIGURE 2.3.

For part vi), to see that $\Pi$ is closed under inverses, suppose $w \in \Pi$, non-zero, $w = re^{it}$. Since $r^2$, which is the sum of the squares of the coordinates belongs to $X = Y$, by i) and ii), $1/r^2$ is constructible, by Lemma 2.7. Reflecting $w$ across the x-axis by Corollary 2.4 gives $re^{-it}$; since, $1/r^2$ can be constructed and so also, by Corollary 2.3, $re^{-it}/r^2$ which is $1/w$. If $w$ and $z$ are constructed then to construct $wz$, we use the simple observation that $2wz = (w+z)^2 - w^2 - z^2$; thus it suffices to construct $w^2$ for any constructible $w$. In order to do that we reflect the complex number 1 about the line from the origin to $w$. This gives us a point $u$ of unit length on the same line as $w^2$, so that by multiplying (the coordinates) by $r^2$ we get $r^2 u = w^2$.

If a line through the origin is constructed with slope $\mu$ then the intersection with the vertical line through $(1, 0)$ gives the point with $\mu \in Y$; the converse follows, since $(1, \mu)$ is constructible. Thus the set of elements in $Y \cup \infty$ is just the set of constructible slopes. ■

For a given $z = a+bi$, the bisectors of constructed segments belong to the field generated by $a, b$ and $i$ over $\mathbb{Q}$. This is easily shown inductively, for points constructed initially from $0, 1, z$. Also, any constructed line has its slope in $\mathbb{Q}(a, b)$. Hence the coordinates of any constructed point belong to $\mathbb{Q}(a, b)$. Thus for all non-real $z$, $\Pi[z] \subseteq Q(a, b, i)$. Also using Corollary 2.8, we have that for any constructed point $w = (x, y)$, $w^2 = (x^2 - y^2, 2xy)$ is also constructible. Since we can invert the squared length, $r^2 = x^2 + y^2 \in X$ of a constructed point $(x, y)$, by Corollary 2.8, the next Corollary follows easily from previous remarks.

**Corollary 2.9.** *For any non-real complex $z=a+bi$, $\Pi[z]$ is a field over $\mathbb{Q}$ containing $z$, closed under complex conjugation, and contained in the field $\mathbb{Q}(a, b, i)$.*

**Corollary 2.10.** *Suppose that $z=a+bi$, with $b$ (real) algebraic. If the irreducible integer polynomial $q(t)$ satisfied by $b$ is not a polynomial in $t^2$ then $X=Y$. Hence, if $b$ has odd degree or $z$ is a non-real algebraic number of odd degree then $X=Y$.*

**Proof.** Collecting together the even terms gives a non-zero integer polynomial expression in even powers of $b$, which is in $X$. This is equal to a integer polynomial expression in odd powers of $b$ which is in $Y$. Thus $X$ and $Y$ meet so we get the field $X(i)$ by Corollary 2.8. If $z$ has odd degree, the field $\mathbb{Q}(z, \bar{z})$ also has odd degree while $\mathbb{Q}(z, \bar{z}, i) =$ $\mathbb{Q}(a, b, i)$ has even degree, so $\mathbb{Q}(z, \overline{z})$ is not equal to $\mathbb{Q}(z, \overline{z}, i)$. Hence, since $b$ is real, its degree is odd and thus $X = Y$, by the remarks above. ∎

2.2. **Thalian Numbers.** For any non-real complex number $z = a+bi$, the field $X$ contains $\mathbb{Q}(a, b^2)$ and is contained in $\mathbb{Q}(a, b)$; also, $Y = Xb$.

The non-real complex number $z$ is termed Thalian if $\Pi[z]$ contains $i$. Hence, for $z$ Thalian, $\Pi[z]$ is the field $X(i)$, where $X = Y$ and conversely. For any Thalian, since $b \in X$, then $X = \mathbb{Q}(a, b)$ and $\Pi[z] = \mathbb{Q}(a, b, i) = X(i)$. For a non-Thalian complex number, $X$ contains $\mathbb{Q}(a, b^2)$ and since $\Pi[z]$ does not contain $i$, it follows that $b \in Y$ and not in $X$; hence, $Y = Xb$. Consequently, $\Pi[z] = \mathbb{Q}(a, b^2) + \mathbb{Q}(a, b^2)bi$. Furthermore, $\Pi[z]$ is the field $\mathbb{Q}(z, \overline{z})$.

Summarizing these results, we have the following.

**Theorem 2.11.** *For any non-real complex $z=a+bi$, $\Pi[z]$ is a field of degree 2, generated by $bi$, over the field $X$. The field $X = \Pi[z] \cap \mathbb{R}$ is either $\mathbb{Q}(a, b^2)$ or $\mathbb{Q}(a, b)$. The non-real complex number is Thalian iff $b \in \mathbb{Q}(a, b^2)$.*

**Corollary 2.12.** *The field $\mathbb{Q}(z, \overline{z})$ belongs to $\Pi[z]$. If $z$ is non-Thalian, $\Pi[z]$ is the smallest field containing $z$ which is also closed under complex conjugation. If $z$ is Thalian, $\Pi[z]$ is the smallest field containing $z$ and $i$, which is also closed under conjugation and conversely.*

**Proof.** If $z$ is non-Thalian, $\mathbb{Q}(a+bi, a-bi) = \mathbb{Q}(a, bi) = \mathbb{Q}(a, b^2, bi) = \Pi[z]$. If z is Thalian and then $b \in \mathbb{Q}(a, b^2)$, so $\mathbb{Q}(a, b^2, bi) = \mathbb{Q}(a, b, i) = \mathbb{Q}(a+bi, a-bi, i)$. It may happen that a field containing $z$ and closed under conjugation, also contains $i$; in that case, $z$ is Thalian. ∎

**Corollary 2.13.** *If $D < 1$ is a square-free integer, $z=\sqrt{D}$ is a non-Thalian.*

**Corollary 2.14.** *For $m>2$, the root of unity $e^{i2\pi/m}$ is a Thalian iff $m$ is divisible by 4.*

**Proof.** Let $z = e^{i2\pi/m}$. If $m$ is divisible by 4 then $z^{m/4} = i$ and thus $z$ is a Thalian. For the converse, we consider the case of $m$ odd or twice an odd integer. For any integer m>1, $\mathbb{Q}(a, b^2) = \mathbb{Q}(cos(2\pi/m), sin(2\pi/m)^2) = \mathbb{Q}(cos(2\pi/m))$ has degree $\phi(m)/2$; whereas, for an integer not divisble by 4, $\mathbb{Q}(b) = \mathbb{Q}(sin(2\pi/m))$ has degree $\phi(m)$, by Lehmer [11]. Therefore $\mathbb{Q}(a, b) = \mathbb{Q}(cos(2\pi/m), sin(2\pi/m))$ has degree at least $\phi(m)$, so is different from $\mathbb{Q}(a, b^2)$. Hence $z$ is a non-Thalian. ∎

## 3. Pythagorean Constructions and Numbers

In order to get more structure on $\Pi$, we add another axiom:

(4) the line bisecting any given constructed angle can be constructed. The axioms (0)-(4) are the axioms for Pythagorean constructions. We shall assume that the first two given points in axiom (0) are 0 and 1.

**Theorem 3.1.** *Given the axioms (0)-(3), the following are equivalent:*
  *i) a unit length segment can be marked on any constructed ray;*
  *ii) the angle bisector axiom (4);*
  *iii) a constructed segment's length can be marked on any other constructed ray.*

**Proof.** iii)$\to$ i) is trivial; i)$\to$ ii) We first show that all constructed angles can be bisected, if we can mark unit length. Without loss of generality, we can assume the angle is (strictly) less than 180 degrees. As before we can move the angle to the origin and mark unit lengths on each ray of the angle, giving us points $A$ and $B$. Next drop perpendiculars at $A$ and $B$, which will meet at a point $C$ (since the angle is less than 180 degrees) The 4 points $A, B, C$ and the origin 0 give us two right triangles with a common hypotenuse and two legs meeting at 0 of length 1. Therefore the third sides are equal and so the triangles are congruent; thus, the angle has been bisected. ii)$\to$ iii) Suppose that all constructed angles can be bisected. Suppose that a segment length $AB$ has been constructed and a ray $L$ starting at $C$ is given. We can first move segment $AB$ so that it starts at $C$, say $CD$ by Lemma 2.1. Now bisect the angle between $CD$ and the given ray $L$. Next reflect $CD$ across this bisector line, so that we have now marked the length $AB$ along the given ray. ∎

**Corollary 3.2.** *The angle at the origin formed by 1, 0 and $z = e^{i\pi/k}$, $k$ an odd integer $\neq 1$, can be bisected by Thalian constructions in $\Pi[z]$, but a unit length on that ray is not in $\Pi[z]$.*

**Proof.** From the discussion above we can clearly bisect the angle, since $z$ and 1 have length 1. However if a unit length could be marked then $\Pi[z]$ would contain a $4k$th root of unity, and hence $i$; however $z$ is a non-Thalian so $\Pi[z]$ does not contain $i$; and consequently also, this constructed bisector has a length which does not belong to $X$. ∎

One can then view this axiom as providing the 'tool' to construct the points on a circle centered at a constructible point $P$, with a given constructible radius $r$, by constructing points on the lines through $P$ with any constructible slope. We denote by $\pi$, the constructed numbers using axioms (0)-(4). With axiom (0) we use starting points $0, 1, i$. Actually $i$ is constructed from the axioms, (0)-(4) just starting with 0 and 1, since the bisection of the 90 degree angle between the axes then gives by reflection the unit direction on the y-axis. Also, this shows that $X = Y$. The other important consequence (of axiom (4)) is that lengths of segments can be constructed; that is the field $X$ is closed under $\sqrt{a^2 + b^2}$ for every constructed point $(a, b)$. Thus $X$ is the Pythagorean numbers, $\mathcal{P}$, the smallest subfield of the real numbers containing $\mathbb{Q}$ which is closed under the operation of taking square roots of a sum of two squares, i.e., closed under $\sqrt{1 + x^2}$. Consequently, we have the following.

**Theorem 3.3.** *The constructible points in $\pi$ is the field $\mathcal{P} \oplus i \mathcal{P}$. The positive elements of $\mathcal{P}$ is the set of constructible segment lengths.*

In fact as Hilbert shows, [6], Theorem 66, the collection of all (positive) real elements of this Pythagorean field is the same as the set of all real algebraic numbers which are obtained by extraction of (positive) square roots, and which are totally real, meaning that any of its algebraic conjugates are also real numbers. It is easy to see that the positive Pythagorean reals are totally real since they are obtained from $\mathbb{Q}$, using only field operations, and the operation $\sqrt{1 + x^2}$, which is a totally real number for any totally real number $x$. The converse, that the totally real numbers which have degree a power of 2 are real Pythagoreans is more complicated; the details can be found in [1].

The number $\sqrt{1 + (1 + \sqrt{2})^2} = \sqrt{4 + 2\sqrt{2}}$ is a Pythagorean number since $\sqrt{2}$ is Pythagorean. However, a right triangle with a (horizonatal) leg of length $\sqrt{2}$ and a hypotenuse of length $\sqrt{4 + 2\sqrt{2}}$ would have a (vertical) leg of length , $\sqrt{2 + 2\sqrt{2}}$; it is not totally real since its conjugate $\sqrt{2 - 2\sqrt{2}}$ is not real. Thus the field of real Pythagorean numbers is closed under $\sqrt{a^2 + b^2}$, but not $\sqrt{a^2 - b^2}$, even when $a^2 > b^2$.

## 4. Euclidean Constructions and Numbers

We next add the axiom:
(5) Given a constructed line $l$ and constructed points $P, Q$, then whenever possible, the line through $Q$, which reflects $P$ onto $l$, can be constructed.

One can achieve this origami construction, by folding $P$ to $l$, and then sliding $P$ along $l$, until the fold line passes through $Q$.

The axioms (1)-(5) are the axioms for the Euclidean origami numbers, $\mathcal{E}$, obtained using the axioms beginning with the set containing just 0 and 1. These Euclidean origami constructions (1)-(5) enable us to construct exactly the same set as we could obtain by using the ruler and compass constructions. As is well known, the field of numbers constructed by ruler and compass is the smallest field containing $\mathbb{Q}$ and closed under taking square roots. The process of taking square roots of the complex number $z = re^{i\theta}$, can be viewed as involving the two steps of bisection of the angle $\theta$ and the square root of $r$. Certainly we can bisect any constructed angle just using axiom (4). Also, we can take extract some square roots of numbers from the Pythagorean field, but not all.

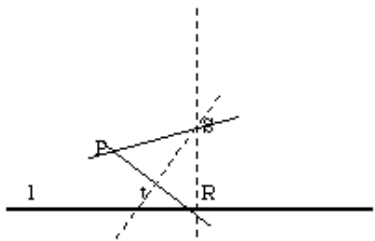

FIGURE 4.1.

We next develop some consequences of axiom (5), and thereby obtain closure under square roots. Consider a parabola, $\mathcal{K}$, having directrix $l$ and focus $P$. Axiom 5 allows us to construct the points of this parabola and the tangent lines there. To see this, use axiom 5 with the focus $P$ and directrix $l$ to construct the line $t$ passing through some auxiliary point $Q$ which reflects $P$ onto $l$. The perpendicular to $t$ passing through $P$ meets the given line $l$ in this constructible point $R$. Next, construct the perpendicular to $l$ at $R$. The intersection of this perpendicular, and the line $t$ is the point $S$.

Since $t$ is the perpendicular bisector of $PR$, the point $S$ is equidistant from $P$ and $R$, so is on the parabola $\mathcal{K}$; the line $t$ is tangent to the $\mathcal{K}$ at S, since it bisects the angle $PSR$ and therefore satisfies the equal angles characterization property, that a line $t$ is a tangent of a parabola, $\mathcal{K}$, if $t$ bisects the angle formed by the lines $PS$ and the line parallel to the axis of $\mathcal{K}$ through $S$.

It is now easy to use a parabola to construct square roots. Let $P = (0, 1)$, and use the directrix $l$, $y = -1$, then the parabola has the equation $y = \frac{1}{4}x^2$. The tangent line to this parabola at the point, $(x_0, \frac{1}{4}x_0^2)$, has slope $m = \frac{1}{2}x_0$; the tangent line has equation $y - \frac{1}{4}x_0^2 = \frac{1}{2}x_0(x - x_0)$. The intersection of this tangent with the line $x = 0$ gives the point $Q = (0, -\frac{1}{4}x_0^2)$. Therefore we use axiom 5, with the focus $P = (0, 1)$, directrix $y = -1$, and auxiliary point $Q = (0, \frac{-1}{4}r)$, to construct the point on the parabola having x-coordinate $\sqrt{r}$.

Consequently, we can construct all square roots of complex numbers in the field $\mathcal{E}$; also, it is easy to see that any new point constructed by using axiom 5 uses only field operations and square roots of previously constructed numbers. Thus the field obtained from using these origami axioms is just precisely the Euclidean constructible complex numbers.

**Theorem 4.1.** *The constructible points in $\mathbb{C}$, obtained by using axioms (1)-(5), starting with the numbers 0 and 1, the field of Euclidean constructible numbers is the smallest subfield of $\mathbb{C}$, closed under square roots.*

It is easy to see that the real subfield of $\mathcal{E}$ is the smallest real subfield closed under taking square roots of its positive elements. As a final comment here we show axioms (2), (3) and (5) are equivalent to axioms (1)-(5). First, to deduce axiom (4), given constructed lines $l$ and $m$, we can construct the point of intersection $Q$ and choose another constructible point $P$ on $l$. The lines constructed by axiom (5) which pass though $Q$ and reflect $P$ onto $m$ will bisect the angles at $Q$. To deduce axiom (1), we are given two constructible points $P, Q$ and by axiom (3) we can construct the perpendicular bisector $l$ of the segment $PQ$. By (5) we can construct a line $m$ through $Q$ which reflects $P$ onto $l$. Next by (5) we can construct a line through $P$ which reflects $Q$ onto $m$; but since $Q$ is already on $m$, this is the line through $P$ and $Q$. Summarizing these comments we obtain the following result.

**Theorem 4.2.** *The field of Euclidean constructible numbers $\mathcal{E}$ is the smallest subset of constructible points in the complex numbers $\mathbb{C}$, which contains the numbers 0 and 1, and is closed under the axioms:*

*($\alpha$) The point of coincidence of two constructible lines is a constructible point.*

*($\beta$) The perpendicular bisector of the segment connecting two constructible points is a constructible line.*

*($\gamma$) Given a constructed line $l$ and constructed points $P, Q$, then whenever possible, the line through $Q$, which reflects $P$ onto $l$, can be constructed.*

## 5. Conic Constructions and Origami Numbers

We add the final axiom:

(6) Given constructed lines $l, m$ and constructed points $P, Q$, then whenever possible, any line which simultaneously reflects $P$ onto $l$ and $Q$ onto $m$, can be constructed.

A MATHEMATICAL THEORY OF ORIGAMI CONSTRUCTIONS AND NUMBERS

One can achieve this origami construction, by folding $P$ to $l$, and then if possible, sliding $Q$ until it lies along $m$. This construction is the simultaneous tangent line to the two parabolas described with the given data of directrices and foci. It certainly can not always be accomplished, since such a line does not necessarily exist.

The axioms (1)-(6) are the (origami) construction axioms for the complex origami numbers, $\mathcal{O}$. The origami constructions (1)-(6) enable us to construct a real solution to a cubic equation with real coefficients in this field $\mathcal{O}$. To see this, consider the conics

$$(y - \tfrac{1}{2}a)^2 = 2bx, \quad y = \tfrac{1}{2}x^2.$$

These conics have foci and directrices that are constructible using field operations involving $a$ and $b$. Consider a simultaneous tangent, a line with slope $\mu$ meeting these curves at the respective points $(x_0, y_0)$, $(x_1, y_1)$. It is important to realize that by Corollary 2.8, a line is constructible iff its slope is a number of the field or $\infty$. Now, differentiation yields, $\frac{b}{y_0 - \frac{a}{2}} = \mu = x_1$, and so $y_1 = \tfrac{1}{2}\mu^2$, $x_0 = \frac{(y_0 - \frac{a}{2})^2}{2b} = \tfrac{1}{2}b\mu^2$; thus $\mu = \frac{y_1 - y_0}{x_1 - x_0} = \frac{\frac{\mu^2}{2} - \frac{a}{2} - \frac{b}{\mu}}{\mu - \frac{b}{2}\mu^2}$. Simplifying we get that $\mu$ satisfies $\mu^3 + a\mu + b = 0$, and hence we can solve any cubic equation with specified real constructible $a, b \in \mathcal{O}$ for its real roots.

For example, the Delian problem was solved this way by Menaechmus, tutor to Alexander the Great, also famous for instructing his pupil '...there are royal roads and roads for commoners, but in geometry there is one road for all'. The Delian problem, or duplication of the cube, requires the construction of $2^{\frac{1}{3}}$, and this is done simply with the parabolas described above, solving $\mu^3 - 2 = 0$.

Furthermore, this enables us to show that we can trisect constructed angles, by solving the Chebychev equation $4x^3 - 3x = cos(3\theta)$ for $x = cos(\theta)$. Pappus had an alternative method using conics ([13]); a pure origami construction for trisection of angles is given in [7].

For example, to construct a regular 9-gon, we must solve the Chebychev equation when $\theta = 40$ degrees; we can solve this equation $4x^3 - 3x = \frac{-1}{2}$, by generating the common tangent to the parabolas $y = \tfrac{1}{2}x^2$ and $(y + \tfrac{3}{8})^2 = \tfrac{1}{4}x$. There are in fact three common tangents as displayed below. It is easy to see, using first and second derivative tests that the real equation $\mu^3 + a\mu + b = 0$ will have three real roots if $27b^2 + 4a^3 < 0$. In the example of the parabolas used for the 9-gon, pictured below on the left, the fourth tangent is the line at $\infty$ in the projective plane containing this ordinary real affine plane. (It is a consequence of Bezout's Theorem for algebraic curves that conics have at most 4 common tangents.)

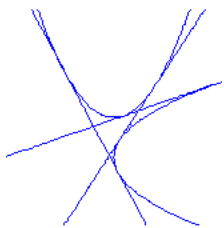

Figure 5.1. 9-Gon Conics

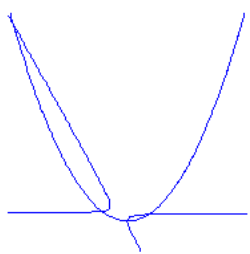

Figure 5.2. 9-Gon Duals

The roots of the general complex cubic with constructible coefficients can also be constructed; one can see this from the explicit solutions (Cardano's formula), to the cubic which involve only square and cube roots. Since we can bisect and trisect constructed angles, and take real square and cube roots with the aid of axiom (6) the roots of a complex cubic with coefficients in $\mathcal{O}$ also belong to $\mathcal{O}$.

In this next section we introduce some classical concepts from algebraic geometry.

5.1. **Higher Geometry.** We extend the plane to the projective plane, by adding in the line at infinity. The points of this projective plane can be viewed as triples of real numbers $(x, y, z)$, not all zero, modulo the scalar multiplications by non-zero reals, i. e., $(x, y, z) \equiv (sx, sy, sz)$, for any non-zero scalar s. The ordinary plane is viewed as the (equivalence classes of) points where $z$ is non-zero. The line at infinity is the equivalence class of points where $z = 0$. An equation which is homogeneous of degree 2 in three variables is called a conic equation, since its solutions describe the points of a conic in the projective plane. One can view this equation in matrix-vector form as

$$F(x, y, z) = (x, y, z)A(x, y, z)^t = 0,$$

where $A$ is a non-zero three by three real symmetric matrix. So far, the origami axiom (6) gives us the construction of simulataneous tangents to two parabolas.

We can use the method of dual curves to allow us to reformulate the simultaneous tangent of two point conics as the common point of the dual line conics. In this way we can use axiom (6) to construct the common intersection points of two parabolas. The payoff of developing some rudimentary projective geometry here is that we will also be able to construct the common points (or common tangents) to two conics. This will lead to the characterization of origami numbers. The dual to the point conic $F(x, y, z) = 0$ is the equation which is satisfied by all the tangents to $F$; viewing the tangent line as $ux + vy + wz = 0$ the coordinates $(u, v, w)$ satisfy the equation $H(u, v, w) = (u, v, w)Adj(A)(u, v, w)^t$, where $Adj(A)$ is the adjoint of $A$, ([5]). This is the line conic. The dual of the line conic is the original point conic. For example, the dual conics to $y = \frac{1}{2}x^2$ and $(y + \frac{3}{8})^2 = \frac{x}{4}$ are $v = \frac{1}{2}u^2$ and $-6uv + 16u = v^2$. These are displayed in the figure above. Corresponding to the point of intersection $(a, b, 1)$ in the projective plane on the right is the linear equation $ax + by + 1 = 0$ in the plane on the left; the origin $(0, 0, 1)$ corresponds to the line at infinity $z = 0$. Thus the picture displays all four intersection points of the parabola and hyperbola.

To get the proper perspective on point conics, one can view the non-degenerate conics in the affine $(x, y, z)$ plane $z = 1$ but should consider also the behavior at infinity, $z = 0$. The parabola has tangent line $z = 0$. The hyperbola has two points (corresponding to the ends of the asymptotes) on $z = 0$. The ellipse does not have a (real) point on $z = 0$. For example the two conics $y = x^2$ and $(x - 1)^2 + (y - 1)^2 = 1$ have a common tangent, the line $y = 0$. We can move this tangent line to $z = 0$, by a linear change of variables. The conics are projectively $yz = x^2$ and $x^2 - 2xz + y^2 - 2yz + z^2 = 0$. After the linear change of variables, permuting $y$ and $z$, we obtain the conics $yz = x^2$ and $x^2 - 2xy + z^2 - 2yz + y^2$, which in the affine plane $z = 1$, gives $y = x^2$ and $x^2 - 2xy + 1 - 2y + y^2$, both parabolas, which have $z = 0$ as a common tangent.

The last topic in classical geometry that I need to bring in, is the notion of a pencil of conics. We consider two real symmetric matrices $A$ and $B$. The pencil of these is the family of real symmetric matrices $A - \lambda B$. If $A, B$ give conic equations $F, G$, respectively, then the equation for a conic in the pencil is $K(x, y, z) = F(x, y, z) - \lambda G(x, y, z) = (x, y, z)(A - \lambda B)(x, y, z)^t$. One of the most important things about pencils is that any conic in the pencil contains the simultaneous solutions to the (or any independent) two generating conics. That is, if $F(x_0, y_0, z_0) = G(x_0, y_0, z_0) = 0$ then also $K(x_0, y_0, z_0) = 0$.

One can solve fourth degree equations with these kinds of methods. Consider the real quartic

$$x^4 + ax^2 + bx + c = 0.$$

Let $y = x^2$. Naively, then we can write the quartic as the simultaneous solution of the two parabolas, when $b \neq 0$,

$$y = x^2, \ (y + \frac{1}{2}a)^2 = -b(x + 4c - \frac{a^2}{4b}).$$

The problem is to find the simultaneous solutions or common points of these parabolas. The common points are the duals of the common tangents of the dual curves. This method involves the cubic which is the determinant of the pencil determined by the dual curves; this is called the resolvant cubic equation of the quartic. The solutions of the resolvant cubic give the conics of the pencil which are degenerate, that is, the conic factors into linear factors (possibly over the complex numbers).

**Lemma 5.1.** *The intersection points and common tangents of two distinct non-degenerate conics with equations defined over $\mathcal{O}$ can be constructed by origami methods defined by axioms (1)-(6).*

**Proof.** We consider the common tangents, if any, to two non-degenerate conics $F, G$ If there is a tangent, then this tangent is also tangent to every conic in the pencil that $F, G$ generate, so it is a common tangent on the degenerate conics in the pencil, when $det(A - \lambda B) = 0$. Solving for $\lambda$ involves solving the cubic determinant equation; but this can be done by origami. By a change of coordinate system, which uses only the field operations, we move this tangent line (by a linear transformation) to the line at infinity in the projective plane. Now we have two conics with a common tangent at infinity; thus, these new conics are parabolas. Any further common tangents now in the new affine plane can be constructed using origami.

We can reduce the problem of construction of common points of two curves to the construction of common tangents for the dual curves. Dualizing the line conic gives the original point conic. The equation for the dual curve is based on the adjoint so its coefficients are also in $\mathcal{O}$. ∎

Consider the set of points in the complex numbers which are obtained as intersections of lines or conics with coefficients in the real subfield $\mathcal{O}_R$ of the origami complex numbers $\mathcal{O}$, those points which are constructible by axioms using (1)-(6). These will be called the *conic constructible* points. This is equivalent (using Theorem 5.2) to the notion of conic constructible points developed in [13]. In that article, constructibility of directrices, eccentricity, foci, radius, etc. are the conditions for conic constructibility.

Summarizing the consequences of this section, we obtain the following theorem.

**Theorem 5.2.** *The constructible points in $\mathbb{C}$, obtained by using axioms (1)-(6), starting with the numbers 0 and 1, is the field of origami constructible numbers, $\mathcal{O}$; it is the smallest subfield of $\mathbb{C}$, closed under square roots, cube roots and complex conjugation. This field $\mathcal{O} = \mathcal{O}_R \oplus i\,\mathcal{O}_R$ is also the set of conic constructible points. The field $\mathcal{O}_R$ is the smallest subfield of the reals closed under arbitrary (real) cube roots and square roots of its positive elements.*

**Proof.** We have seen already that the fields $\mathcal{O}$ and $\mathcal{O}_R$ are closed under conjugation, square and cube roots, either complex or real, respectively. Thus, the smallest such subfields of the complex numbers or real numbers closed under these specified roots and conjugation, $\mathcal{M}$ and $\mathcal{M}_R$, are contained in $\mathcal{O}$ and $\mathcal{O}_R$, respectively.

The conic constructions of common points of distinct conics can all be done using origami constructions by using Lemma 5.1, so the conic constructible points are contained in $\mathcal{O}$ and their coordinates are in $\mathcal{O}_R$. Furthermore it follows immediately from the argument in Lemma 5.1, that these points and their coordinates are obtained by using field operations and solving cubic equations with coefficients in $\mathcal{O}_R$. The coordinates for intersections of lines or line and conic also involve solving either a linear equation or a cubic, possibly reducible, with coefficients in $\mathcal{O}_R$. Hence coordinates of conic constructible points are contained in $\mathcal{M}_R$.

On the other hand, an origami constructible point in $\mathcal{O}$, has its coordinates in $\mathcal{O}_R$, and is easily obtained as the intersection of a horizontal and vertical line with coefficients in $\mathcal{O}_R$. Thus $\mathcal{O}$ is contained in the conic constructible points. So $\mathcal{O} = \mathcal{M}$ is the set of conic constructible points and also $\mathcal{O}_R = \mathcal{M}_R$. Now, since the conic constructible points form a field closed under complex conjugation, the coordinates of any of its points are in $\mathcal{O}_R$ ∎

As a final comment, we notice that the axiom system for $\mathcal{O}$ can easily be simplified using the remarks preceeding Theorem 4.2 to obtain the following result ( [2]).

**Theorem 5.3.** *The field of origami constructible numbers $\mathcal{O}$ is the smallest subset of constructible points in the complex numbers $\mathbb{C}$, which contains the numbers 0 and 1, and is closed under the axioms:*

*($\alpha$) The point of coincidence of two constructible lines is a constructible point.*

*($\beta$ ) The perpendicular bisector of the segment connecting two constructible points is a constructible line.*

*($\gamma$) Given a constructed line $l$ and constructed points $P, Q$, then whenever possible, the line through $Q$, which reflects $P$ onto $l$, can be constructed.*

*($\delta$) Given constructed lines $l, m$ and constructed points $P, Q$, then whenever possible, any line which simultaneously reflects $P$ onto $l$ and $Q$ onto $m$, can be constructed.*

E-mail: alperin@mathcs.sjsu.edu, Department of Mathematics and Computer Science, San Jose State University, San Jose, CA 95192 USA